\documentclass{article}
\usepackage{amsmath,amssymb,amsfonts,amsthm,color,graphicx,epstopdf}
\usepackage{float}

\newtheorem{theorem}{Theorem}
\newtheorem{lemma}[theorem]{Lemma}
\newtheorem{proposition}[theorem]{Proposition}

\DeclareMathOperator{\Mod}{mod}

\begin{document}

\title{\bfseries Cournotian dynamics of spatially distributed renewable resources}

\author{Sebastian Ani\c{t}a$^{a,b}$  \quad Stefan Behringer$^{c,\star }$ \\
        Ana-Maria Mo\c{s}neagu$^{a}$ \quad Thorsten Upmann$^{d}$}

\renewcommand{\thefootnote}{\fnsymbol{footnote}}
\footnotetext{$^a$ Alexandru Ioan Cuza University of Ia\c si, Romania,
        $^b$ Octav Mayer Institute of Mathematics, Ia\c si, Romania,
        $^c$ Sciences Po, Paris, France,
        $^d$ Bielefeld University \& CESifo Munich, Germany.} %
\renewcommand{\thefootnote}{\fnsymbol{footnote}{$\star$}}
\footnotetext{Corresponding author: Sciences Po, Department of Economics,
                        28 rue des Saints-Per\'{e}s, 75007 Paris, France}
                        \renewcommand{\thefootnote}{\fnsymbol{footnote}}
\footnotetext{E-mail adresses:
        sanita@uaic.ro (S.~Ani\c{t}a),
        Stefan.Behringer@sciencespo.fr (S.~Beh\-ring\-er),
        anamaria.mosneagu@uaic.ro (A.~M.~Mo\c{s}neagu),
        tupmann@wiwi.uni-bielefeld.de (T.~Upmann).}

\date{\normalsize\today}

\maketitle

\begin{abstract}
We extend modern Walrasian economics, and in particular the results on
Cournot convergence and dynamics, by focusing on renewable resources in a
spatial setting. Building on the harvesting model of Behringer and Upmann
(2014) we endogenize prices and investigate the two cases of durable and
non-durable renewable commodities. We find that endogenizing prices is
sufficient to prevent the full exploitation result and look at how
competition affects not only the stock but also the temporal incentives for
exploitation. We derive convergence results in static and dynamic settings
which suggest that the classical Cournotian outcomes may prevail.
\end{abstract}




\renewcommand{\thefootnote}{\arabic{footnote}}

\section{Introduction}

The theory of perfect competition originating in the works of Cournot and
Edgeworth has been successfully extended to non-cooperative settings that
have a dynamic nature. Important landmarks in this direction, in particular
by Green and Radner, are assembled in a special edition of the Journal of
Economic Theory 1980. The former has shown that in a \emph{repeated game}\
setting, where a stage game is replicated, a small degree of noise
(imperfect information about some aggregate statistic) is sufficient to get
back to the stationary Cournot outcome as individual deviations from
collusive arrangements cannot be detected with sufficient accuracy. This
\textquotedblleft limit principle\textquotedblright \ holds even if we are
dealing with finitely many agents only. Subsequently, Levine and Pesendorfer
(1995) show that for the collusive outcome to be sustained, the aggregate
noise level has to decrease with the number of agents sufficiently fast, see
also Al-Najjar and Smorodinsky (2000).

While classical mircoeconomic theory deals with homogenous consumption goods
in a static framework, we consider a dynamic framework. More precisely, we
consider the harvesting and sale of a renewable natural resources (fish,
timber, game) the stock of which obeys a given law of growth.\footnote{%
See Smith (1968, 1977), Beddington et.~al.~(1975) or Clark et.~al.~(1979)
for early economic analyses.} In addition, we allow for the resources to be
spatially extended taking into consideration demands from the discipline and
policy makers, see Deacon et al. (1998). We then investigate the validity of
the classical Courotian results for goods that are heterogenous, thereby
extending the classical results to renewable commodities in a dynamic
setting. To this end we endogenize prices for both types of commodities,
non-durable and durable renewable ones by taking into account output market
behaviour.

Recently Behringer and Upmann (2014) investigate optimal harvesting of a
renewable resource that is spatially distributed over a \emph{continuous}
domain. Since the agent is required to move in space, an optimal policy
consists of an optimal choice of both, harvesting and movement. This
approach contrasts with previous analyses of \emph{discrete} spaces, e.g.
Sanchirico and Wilen (1999, 2005) but is similar to Belyakov and Veliov
(2014) who also consider a continuous setting.\footnote{%
Harvesting models have also been intensively studies by Sebastian Ani\c{t}a,
see Ani\c{t}a (2000) and Ani\c{t}a, Ani\c{t}a, and Arn\u{a}utu (2009) and
also the references therein.}

The dynamic optimization problem in the model of Behringer and Upmann (2014)
consists of a simultaneous choice of the speed of movement $\left \{
v(t)\right \} _{t\in \mathcal{T}}$ and the harvesting rate $\left \{
h(t)\right \} _{t\in \mathcal{T}}$. More precisely, the harvesting agent
moves on the periphery of a unit circle on which the resource, with stock $%
f(\cdot )$, is growing according to some growth function $g(\cdot ).$ The
agent's location $s$ is therefore on $\mathcal{S}=[0,2\pi ]$. $\mathcal{T}$
denotes the harvesting period or season $[0,T]$ and harvesting comes at a
cost $C(\cdot )$ that may depend on the speed of the agent and the
harvesting rate.

As the agent cannot harvest more than the entire resource stock at any
particular location, we have $h(t)\leq \max \left \{ 0,f(t,s(t))\right \} .$
Harvesting takes place only at the actual location of the agent $x=s(t)$ and
implies a downward jump in the stock of the resource $f(\cdot ,x)$ at the
set of arrival times of the agent at that location $x:$ $J(x)=\left \{
t_{1}(x),t_{2}(x),...\right \} .$ Accordingly the law of motion for the
stock is%
\begin{equation}
f_{t}(t,x)=g(f(t,x))\text{ }\forall \text{ }t\in \mathcal{T}\backslash
J(x),x\in S  \label{1a}
\end{equation}

\begin{equation}
f(t^{-},x)-f(t^{+},x)=h(t)\forall \text{ }t\in J(x),x\in S  \label{2a}
\end{equation}%
with constant initial level $f(0,x)=f_{0}(x)$ for all $x\in \mathcal{S}.$

Discounting at a rate $\rho \geq 0,~$the agent's problem is%
\begin{equation*}
\max_{\{v,h\}}\, \int_{0}^{T}e^{-\rho t}\left( h(t)-C(v(t),h(t))\right) dt
\end{equation*}%
s.t.
\begin{equation*}
\begin{aligned} &\dot{s}(t) =v(t), &\forall t\in \mathcal{T} \\ &f_{t}(t,x)
=g(f(t,x)), &\forall t\in \mathcal{T} \backslash J(x), x\in \mathcal{S} \\
&f(t^{-},x)-f(t^{+},x) =h(t), &\forall t\in J(x), x\in \mathcal{S} \\ &h(t)
\in H(t) &\forall t\in \mathcal{T} \\ &f(0,x) =f_{0x}, &\forall x\in
\mathcal{S}\\ &s(0) =0. \end{aligned}
\end{equation*}%
The last line implies that w.l.o.g. we let the agent start at $x=0$ on the
periphery.

For any fixed location equation (\ref{2a}) gives a mapping%
\begin{equation*}
f(t_{i}^{+},x)=G(f(t_{i-1}^{+}(x),x),t_{i}(x)-t_{i-1}(x))-h(t_{i}(x))
\end{equation*}%
where $G$ is the solution of the differential equation between two
consecutive impulses, $G(f,0)=f,$ i.e. we have a problem where time and
space of impulses are related, i.e. not a pure \emph{impulse control problem}
as e.g. Yang (2001).\footnote{%
Note that (\ref{1a}) is \emph{autonomous} and does not depend on
time\thinspace $t$ directly but only via $f(\cdot )$. Hence if we integrate
up (\ref{1a}) over the time of two consecutive rounds $t_{i-1}(x)$ and $%
t_{i}(x)$ we get
\begin{equation*}
f(t_{i}^{+}(x),x)=G(f(t_{i-1}^{+}(x),t_{i}(x)-t_{i-1}(x)))
\end{equation*}%
where we due to the ergodic structure we can now replace the space dimension
by the time difference (time it takes for one round) as time and space are
directly related and it is \emph{either} time \emph{or} space that matters.}

Behringer and Upmann (2014) find that with linear growth and constant speed,
the resource will be fully extinguished by the agent by the end of the
planning horizon. As in the early literature on Walrasian economics, this
work treats prices as exogenous however. In order to fully trace out the
welfare economic consequences of trading renewable natural resource
commodities in the spirit of Cournot we endogenize prices in this paper. In
contrast to the stationary structure of repeated games, our analysis dealing
with renewable commodities allows us to investigate the validity of the
Green's \textquotedblleft limit principle\textquotedblright \ in a truly
dynamic context.

Other recent advances in the wake of Green's work are Al-Najjar and
Smorodinsky (2000) and most recently Kalai and Shmaya (2015a, 2015b),
who relax the original complete information setting and, similar to Jehiel and
Koessler (2008) allow for boundedly rational behaviour as well as mixed
strategies.

\section{Non-durable good analysis}

Consider a fixed location $x\in \mathcal{S}.$ Instead of letting the agent
control the harvest $h(t),$ we assume that the agent controls the harvesting
share $\alpha (t),$ (i.e. uses a fishing net with a given mesh size)\ so
that the harvest amounts to $h(t)=\alpha (t)f(t).$ This is the common
formulation in the resource literature:\ e.g. fish is harvested as a share $%
\alpha (t)$ of the stock and so the yield from fishing is multiplicative in
the stock.

We assume that \emph{costs of harvesting} are linear and normalized to zero,
implying a strictly concave per-unit \emph{net} revenue function of the form
$R(\alpha (t))=\alpha (t)(1-\alpha (t)).$ The per unit profit of the
resource is thus increasing with the share of fish put on the market for low
harvesting shares, attains a maximum at $\alpha =1/2,$ and decreases
afterwards reflecting the fact that the market becomes more and more
saturated as $\alpha $ increases. This specification of the net revenue
corresponds to the existence of a single harvesting agent, who supplies the
market as a \emph{monopolist}. In section 4 we will extend this to an \emph{%
oligopolistic context} where multiple harvesting agents supply the market
and are thus in competition.

We assume that the commodity is \emph{non-durable}, and so cannot be stored
but has to be consumed immediately after purchase. Therefore the quantities
supplied to the market do not accumulate over time. The optimal control
problem is then:
\begin{equation}
\max_{\alpha \in \mathcal{A}}\int_{0}^{T}e^{-\rho t}R(\alpha (t))f^{\alpha
}(t)dt  \label{*}
\end{equation}%
where $\mathcal{A}=\left \{ \alpha \in L^{\infty }(0,T);\,0\leq \alpha
(t)\leq 1\, \,a.e.\right \} $ is the set of admissible controls. As in
Behringer and Upmann (2014) we assume \emph{exponential growth} from here
onwards as this simplifies the presence of the economic discounting factor
substantially.

We denote by $f^{\alpha }(t-)$ the level of the renewable resource at some $%
x=\text{mod}(vt,2\pi )$ \emph{just before} harvesting, so just before the
next supply to the market. Likewise the level of the resource \emph{%
immediately after} harvesting is denoted $f^{\alpha }(t+)$ .$~$The initial
level of the resource at $x$ is $f_{0}(x)$, $\forall x\in \mathcal{S}$ as
motivated above.

We denote for any fixed $\alpha $ and any round $l=\left \{ 0,1,...,k\right
\} $ of the $k$ complete rounds\footnote{%
Note that in order to achieve full coherence between noation and semantics
we refer to the \textquotedblleft $lth$ round\textquotedblright \ as\ the
\textquotedblleft $(l+1)th$ round\textquotedblright ,\ as our counting
variable $l$ starts at zero
which however corresponds to the very first round
that agent undertakes. In addition to a notional extension of the density
function beyond the fixed time horizion (see below) this minor linguistic
inaccuracy makes our analysis and its description much less cumbersome than
it otherwise would be.}
the stock of the resource by $f_{l}^{\alpha
}:[0,\theta )\rightarrow
\mathbb{R}
_{+}$ as a function of the time elapsed since the last arrival (at time $%
l\theta ).~$Note that the time necessary to circle around the periphery once
is $\theta \equiv 2\pi /v,$so that the stock (and thus the density) is a
function of the travelling time $\theta $ (or equivalently of speed $v).$

The travelling time for one complete round on the circle equals the duration
between any two consecutive arrivals times at a (any) location and thus
equals the growth time of the resource between two subsequent harvesting
times. Consequently, the stock of the resource depends on the travelling
time $\theta $ (or speed $v)$ and on the harvesting share $\alpha $
according to the law of motion as motivated above.

Then, using the above definitions we obtain
\begin{equation*}
f^{\alpha }(t+)=(1-\alpha (t))f^{\alpha }(t-)
\end{equation*}%
and because of exponential growth at $r$ it also follows that
\begin{equation}
f^{\alpha }\left( \left( t+\theta \right) -\right) =e^{r\theta }(1-\alpha
(t))f^{\alpha }(t-).  \label{1}
\end{equation}%
Equation (\ref{1}) thus states that the density at time $t+\theta $ just
before harvesting equals the original density at $t$ before harvesting, of
which the harvesting share at $t$ has been deducted and which has since
grown according to the exponential growth rate.

As there are $k$ complete rounds until $T$ we have that%
\begin{equation*}
k\theta \leq T<(k+1)\theta ,\quad \text{ }k\in \mathbb{N}.
\end{equation*}%
For convenience, we extend the time horizon beyond the end of the harvesting
period as
\begin{equation*}
f^{\alpha }(t-)=0\quad \text{on}\quad (T,(k+1)\theta
\end{equation*}%
to allow for $k$ complete rounds of supply and a possibly incomplete round
on the circle with the density after $T$ being zero. This vaporizing stock
after $T\,$\ then notionally extends the time horizon but does not affect
the optimization problem. It only relaxes the effect that the fixed time
horizon has on the possibility to treat only integer rounds.

Now for some round $l\in \left \{ 0,1,2,...k\right \} $ on the circle that
takes place at some time interval $t\in \left[ l\theta ,(l+1)\theta \right]
, $ we define%
\begin{equation*}
f_{l}^{\alpha }\left( \left( t-\theta l\right) -\right) \equiv f^{\alpha
}(t-),\quad l\in \left \{ 0,1,\dots ,k\right \} ,
\end{equation*}%
the stock of the resource just before harvesting extended $l\in \left \{
0,1,2,...k\right \} $ periods into the past. We can then also define the
stock of the resource $l$ periods into the future (by adding time $\theta l$
to the above) as%
\begin{equation*}
f_{l}^{\alpha }(t-)\equiv f^{\alpha }\left( \left( t+\theta l\right)
-\right) .
\end{equation*}%
for any round $l\in \left \{ 0,1,2,...k\right \} $.

Adding time $\theta l$ to (\ref{1}) we find
\begin{eqnarray*}
f^{\alpha }\left( \left( t+\theta +\theta l\right) -\right)
&=&f_{l+1}^{\alpha }\left( \left( t\right) -\right) =e^{r\theta }(1-\alpha
(t+\theta l))f^{\alpha }((t+\theta l)-) \\
&=&e^{r\theta }(1-\alpha (t+\theta l))f_{l}^{\alpha }((t)-)
\end{eqnarray*}%
by the ergodic structure obtained which holds for all $l\in \left \{
0,1,2,...k\right \} $ as $\alpha $ does not impact $f$ differently over
rounds and we make use of the extended time horizon. We thus have for the
time interval $t\in \left[ l\theta ,(l+1)\theta \right] $ that%
\begin{equation*}
f_{l+1}^{\alpha }\left( \left( t\right) -\right) =e^{r\theta }(1-\alpha
(t+\theta l))f_{l}^{\alpha }((t)-)
\end{equation*}%
i.e. the density just before harvesting at any round $l+1$ is given by the
original density in round $l$ just before harvesting, of which the
harvesting share in that round$~$has been deducted and which since has grown
(for one round of time) according to the exponential growth rate. For the
first period, where previous harvesting trivially cannot have a consequence
for present harvest and hence $\alpha $ is not an argument to be
considered,\ this reduces to%
\begin{equation*}
f_{0}^{\alpha }\left( t-\right) =e^{rt}f_{0}(tv)
\end{equation*}%
where $x=$mod$(vt,2\pi )=vt$ if $t\in \lbrack 0,\theta )$ gives the location
in the first round$.$ Thus we find the ergodic relation between round $l\in
\left \{ 0,1,2,...k\right \} $ densities and the following round densities
for $t\in \left[ l\theta ,(l+1)\theta \right] $ as:%
\begin{equation}
\left \{
\begin{array}{c}
f_{l+1}^{\alpha }\left( t-\right) =e^{r\theta }(1-\alpha (t+\theta
l))f_{l}^{\alpha }(t-) \\
f_{0}^{\alpha }\left( t-\right) =e^{rt}f_{0}(tv)%
\end{array}%
\right. .  \label{den}
\end{equation}

The optimal control problem given in \eqref{*} is
\begin{equation}
\max_{\alpha \in \mathcal{A}}G(\alpha )=\max_{\alpha \in \mathcal{A}%
}\int_{0}^{T}e^{-\rho t}\alpha (t)(1-\alpha (t))f^{\alpha }(t)dt.  \label{**}
\end{equation}%
where $\mathcal{A}=\left \{ \alpha \in L^{\infty }(0,T);\,0\leq \alpha
(t)\leq 1\, \,a.e.\right \} $ is the set of admissible controls.

This objective can be rewritten as the sum of $k$ completed and a possibly
incomplete round on the circle as%
\begin{equation*}
\begin{aligned} G(\alpha) =&\sum_{l=0}^{k-1}\int_{0}^{\theta}e^{-\rho
(t+\theta l)}\alpha \left(t+\theta l\right)\left( 1-\alpha \left(t+\theta
l\right)\right) f_{l}^{\alpha }(t-)dt \\ &+\int_{0}^{T-\theta k}e^{-\rho
(t+\theta k)}\alpha \left(t+\theta k\right)\left( 1-\alpha \left(t+\theta
k\right)\right) f_{k}^{\alpha }(t-)dt \end{aligned}
\end{equation*}%
which amounts to a monopoly analysis.

Existence of an optimal control has been shown by Ani\c{t}a, Arn\u{a}utu,
and Capasso (2011), Arn\u{a}utu and Neittaanm\"{a}ki (2003), and Barbu
(1994), going back to a problem of Brokate (1985). Let $\alpha ^{\ast }$ be
such an \emph{optimal} control. Then, for any $w\in L^{\infty }(0,T)$ such
that only $0\leq \alpha ^{\ast }(t)+\varepsilon w(t)\leq 1~a.e.,$ for
sufficiently small $\varepsilon >0$ holds, we have that
\begin{equation*}
G(\alpha ^{\ast })\geq G(\alpha ^{\ast }+\varepsilon w).
\end{equation*}%
Whence we have, making use of the extended time horizon that, summing over
the $k+1$ rounds
\begin{equation}
\begin{aligned} &\sum_{l=0}^{k}\int_{0}^{\theta}e^{-\rho (t+\theta l)}\alpha
^{\ast }\left(t+\theta l\right)\left( 1-\alpha ^{\ast }\left(t+\theta
l\right)\right) f_{l}^{\alpha ^{\ast }}(t-)dt \\
\geq&\sum_{l=0}^{k}\int_{0}^{\theta}e^{-\rho (t+\theta l)}(\alpha ^{\ast
}+\varepsilon w)\left(t+\theta l\right)( 1-\alpha ^{\ast }-\varepsilon
w)\left(t+\theta l\right) f_{l}^{\alpha ^{\ast }+\varepsilon w}(t-)dt.
\label{0} \end{aligned}
\end{equation}%
holds. The following Lemma implies that due to its ergodic structure we can
derive a system \emph{without impulses}.

\begin{lemma}
It holds that
\begin{equation}
\begin{aligned} 0\geq \sum_{l=0}^{k}\int_{0}^{\theta}&e^{-\rho (t+\theta
l)}\left[ w\left(t+\theta l\right)\left( 1-2\alpha ^{\ast }\left(t+\theta
l\right)\right) f_{l}^{\alpha ^{\ast }}(t-)\right.\\ &\left.+\alpha ^{\ast
}\left(t+\theta l\right)\left( 1-\alpha ^{\ast }\left(t+\theta
l\right)\right) z_{l}(t)\right]dt\end{aligned}  \label{3}
\end{equation}%
with%
\begin{equation}
\left \{ \begin{aligned} &z_{l+1}(t) = e^{r\theta}\left[ -w\left(t+\theta
l\right)f_{l}^{\alpha ^{\ast }}(t-)+\left( 1-\alpha ^{\ast }\left(t+\theta
l\right)\right) z_{l}(t)\right], \\ &\qquad \qquad \qquad \qquad \qquad
\qquad \qquad t\in \left[0,\theta \right),\quad l=0, 1, \dots, k-1, \\
&z_{0}(t)=0. \end{aligned}\right.  \label{4.1}
\end{equation}%
where
\begin{equation}
z_{l}=\lim_{\varepsilon \rightarrow 0}\frac{f^{\alpha ^{\ast }+\varepsilon
w}-f^{\alpha ^{\ast }}}{\varepsilon }\quad \text{in}\quad L^{\infty }(0,T).
\label{Diff}
\end{equation}
\end{lemma}

\begin{proof}
See Appendix.
\end{proof}

\subsection{Duality}

We now denote the adjoint state by $p=p(t)$, i.e. $p$ satisfies
\begin{equation}
\left \{ \begin{aligned} &\begin{aligned} &p_{l}(t) =e^{r\theta}\left(
1-\alpha ^{\ast }\left(t+\theta l\right)\right) p_{l+1}(t) \label{5} \\
&\qquad \quad +e^{-\rho (t+\theta l)}\alpha ^{\ast }\left(t+\theta
l\right)\left( 1-\alpha ^{\ast }\left(t+\theta l\right)\right), \\ &\qquad
\qquad t\in \left[0,\theta \right), \quad l = 0, 1, \dots, k-1,
\end{aligned}\\ &p_{k}(t)=\left \{ \begin{array}{l} e^{-\rho (t+\theta
k)}\alpha ^{\ast }\left( t+\theta k\right) \left( 1-\alpha ^{\ast }\left(
t+\theta k\right) \right),\\ \qquad \qquad \qquad \qquad \qquad t\in
\left[0, T-\theta k\right), \\ 0, \quad t\in \left[T-\theta k, \theta
\right]. \end{array}\right. \end{aligned}\right.
\end{equation}%
For the construction of the adjoint problems in optimal control theory we
refer to Ani\c{t}a, Arn\u{a}utu, and Capasso (2011),
Arn\u{a}utu, Neittaanm\"{a}ki (2003), and Barbu (1994).

Defining $a_{l}(t)\equiv \frac{1}{2}\left( 1-e^{\rho \left( t+\theta
l\right) +r\theta }p_{l+1}(t)\right) $, for $t\in \left[ 0,\theta \right) $
and $l=0,1,\dots ,k-1$ we can show the following:

\begin{proposition}
The optimal control $\alpha ^{\ast }$ can be characterized as:%
\begin{equation*}
\alpha ^{\ast }\left( t+\theta l\right) =\left \{
\begin{array}{cc}
a_{l}(t) & \text{for }a_{l}(t)\in \lbrack 0,1] \\
0 & \text{for }a_{l}(t)<0 \\
1 & \text{for }a_{l}(t)>1%
\end{array}%
\right.
\end{equation*}%
for $t\in \left[ 0,\theta \right) $ and $l=0,1,\dots ,k-1$.

In particular, we find for the (potentially incomplete) round $k+1$ that
\begin{equation*}
\alpha ^{\ast }\left( t+\theta k\right) =\frac{1}{2},\quad \text{for}\quad
t\in \left[ 0,T-\theta k\right] .
\end{equation*}
\end{proposition}

\begin{proof}
See Appendix.
\end{proof}

The proof of this proposition makes use of the dual formulation. As argued
above, with the law of motion being the spatial dimension of the problem is
absorbed into the temporal one without loss of generality. The proof now
shows that one may integrate up the system (\ref{4.1}) over consecutive
rounds instead of the whole time horizon and thus segment the problem
further. Then, as system (\ref{4.1}) gives $z_{0}(t)=0$ for the limit of the
differential quotient of the densities for small deviations from the optimum
(\ref{Diff}) for very first round (see footnote 4), one can use the
inequality (\ref{3}) from the above Lemma 1 to find a condition similar to a
\textquotedblleft standard first-order-condition\textquotedblright \ of
optimization despite the complications resulting from the Lebesgue
generalization. The satisfaction of this condition yields the above
proposition.

We thus find that the result of Behringer and Upmann (2014) of full resource
exploitation can be proved more generally. Here however the price effect
will imply that half of the resource is saved. Having a share larger than
half cannot be optimal as it implies that by choosing $1-\alpha $ one may
generate the same revenue $R(\alpha (t))$ but leave more of the resource in
place, which is clearly better. Also the duality analysis allows for
numerical tests that extend the present framework to more realistic and
heterogenous distributions of the resource. This is ongoing work in
Ani\c{t}a et.~al.~(2016).
There it can be seen easily that the standard case for the
non-terminal phases of the horizon is where the adjoint state is
\textquotedblleft large\textquotedblright \ so that $a_{l}(t)<0$ and so it
is optimal for the agent to refrain from harvesting and supplying to the
market. Interestingly as shown in
Ani\c{t}a et.~al.~(2016),
heterogenous
resource distributions may imply that multiple \emph{optimal} control levels
satisfy $a_{l}(t)\in \lbrack 0,1]$ towards the end of the horizon which
suggest that the monopoly harvesting problem becomes more intricate.

\subsection{Aggregate revenue for some constant shares}

By fixing the shares that the monopolist can deliver to the market in each
period (e.g. if there are regulations of minimum mesh size for fisheries) we
can describe the form of the total aggregate revenue in more detail.

For a fixed location $x\in \mathcal{S}$, let $y_{0}=e^{(r-\rho )t}f_{0}(x)$
denote the stock at the first harvest at $x,$ where $f_{0}(x)$ is uniformly
distributed on the periphery. Assuming $N$ rounds, the objective function %
\eqref{**}$~$becomes:
\begin{equation*}
G(\alpha )=\sum_{n=1}^{N}y_{0}\alpha _{n}(1-\alpha _{n})e^{\theta
(n-1)(r-\rho )}\prod_{i=1}^{n-1}(1-\alpha _{i}),
\end{equation*}%
where $\alpha _{n}=\alpha |_{\left[ \theta (n-1),\theta n\right) }$.

\begin{proposition}
 Total revenue with constant $\alpha $ and net growth $\sigma \equiv
r-\rho $ can be calculated as:%
\begin{eqnarray}
&&G(\theta ,\alpha )=\alpha \frac{2\pi y_{0}}{\sigma \theta }\times
\label{Gd} \\
&&\left(
\begin{array}{c}
\frac{1}{-(\alpha +e^{-\sigma \theta }-1)^{2}}\times \\
\left(
\begin{array}{c}
(\alpha -1)(2e^{-\sigma \theta }-e^{-2\sigma \theta }-1)+ \\
\left(
\begin{array}{c}
(\alpha -1+N\alpha )e^{\sigma \theta (N+1)}+ \\
(N\alpha ^{2}-2\alpha +2-2N\alpha )e^{\sigma \theta (N+2)}+ \\
(\alpha -1-N\alpha ^{2}+N\alpha )e^{\sigma \theta (N+3)}%
\end{array}%
\right) e^{-3\sigma \theta }(1-\alpha )^{N}%
\end{array}%
\right) + \\
\left( 1-(N+1)\alpha \right) \alpha ^{N}
\left( e^{\sigma \Mod(T,\theta)}-1\right)
e^{N\sigma \theta }
\end{array}%
\right)  \notag
\end{eqnarray}
\end{proposition}

\begin{proof}
See Appendix.
\end{proof}

This characterization allows for further numerical examples that show the
trade-offs between the speed with which an agent moves on the periphery and
the optimal amounts of the renewable non-durable commodity that is put on
the market.

\subsubsection{Examples}

Given the parameters $T=10,\sigma =3/20,y_{0}=1.$

We first assume that the time to circle once is $\theta =5,$ so that exactly
two rounds are completed as $N=\left \lfloor \frac{T}{\theta }\right \rfloor
=\left \lfloor \frac{10}{5}\right \rfloor =\allowbreak 2.$ W$\text{e have mod%
}(T,\theta )=0,$ so that the final term in (\ref{Gd}) falls out. We then have
\begin{equation*}
\begin{aligned} G(\theta ,\alpha ) =&\sum_{i=1}^{N}E(i)+E(N+1,s(T))=
\frac{8}{3}\pi \alpha \frac{\alpha -1}{\left( \alpha
+e^{-\frac{3}{4}}-1\right) ^{2}}\left(e^{-\frac{3}{2}}
-2e^{-\frac{3}{4}}+2e^{\frac{3}{4}}\right.\\ &-e^{\frac{3}{2}}-8\alpha
e^{\frac{3}{4}} +4\alpha +4\alpha e^{\frac{3}{2}} +8\alpha
^{2}e^{\frac{3}{4}} -\left.2\alpha ^{3}e^{\frac{3}{4}}-3\alpha ^{2}-5\alpha
^{2}e^{\frac{3}{2}}+2\alpha ^{3}e^{\frac{3}{2}}\right) \end{aligned}
\end{equation*}

\begin{figure}[H]
\begin{center}
\includegraphics[scale=0.5]{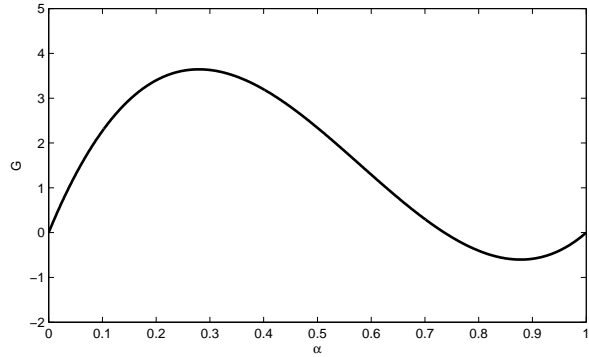}
\caption{$G$ evaluated at $\theta = 5$.}
\label{Fig1}
\end{center}
\end{figure}

We can plot the $G$ function as in Figure~\ref{Fig1} and find the optimal $\alpha$
for two rounds at about $\alpha ^{\ast }=0.28$. Note that, as argued before,
a constant share larger than one half is never optimal and may even yield a
negative objective $G$. With two rounds to complete harvesting everything
already in the first round and thus flooding the market yields a slightly
better outcome than leaving little for the next round but being careful with
the resource by choosing to harvest about a quarter each round yields the
highest objective outcome for the monopoly.

We now assume that the time to circle once is\ $\theta =10,$ so that exactly
one round is completed as $N=\left \lfloor \frac{10}{10}\right \rfloor
=\allowbreak 1,\,$and again mod$(T,\theta )=0$. We have
\begin{equation*}
\begin{aligned} G(\theta ,\alpha ) =&\sum_{i=1}^{N}E(i)+E(N+1,s(T))
=\frac{4}{3}\pi \alpha \frac{\alpha -1}{\left( \alpha
+e^{-\frac{3}{2}}-1\right) ^{2}}\left(e^{-3}-3e^{-\frac{3}{2}}\right.\\
&-e^{\frac{3}{2}}+ 2\alpha e^{-\frac{3}{2}}\left.+2\alpha
e^{\frac{3}{2}}-4\alpha -\alpha ^{2}e^{\frac{3}{2}}+\alpha ^{2}+3\right)
\end{aligned}
\end{equation*}

\begin{figure}[H]
\begin{center}
\includegraphics[scale=0.6]{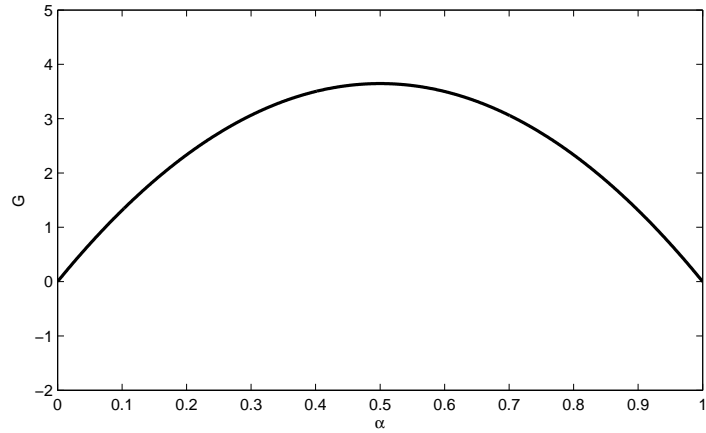}
\caption{$G$ evaluated at $\theta = 10$.}
\label{Fig2}
\end{center}
\end{figure}

We can plot the $G$ function as in
Figure~\ref{Fig2}
and find the optimal $\alpha $
to be $\alpha ^{\ast }=0.5$. Hence with only one round to complete the
monopolist simply supplies the optimal static quantity of one half to the
market.

\subsection{Aggregate revenue for any constant shares}

\begin{figure}[H]
\begin{center}
\includegraphics[scale=1]{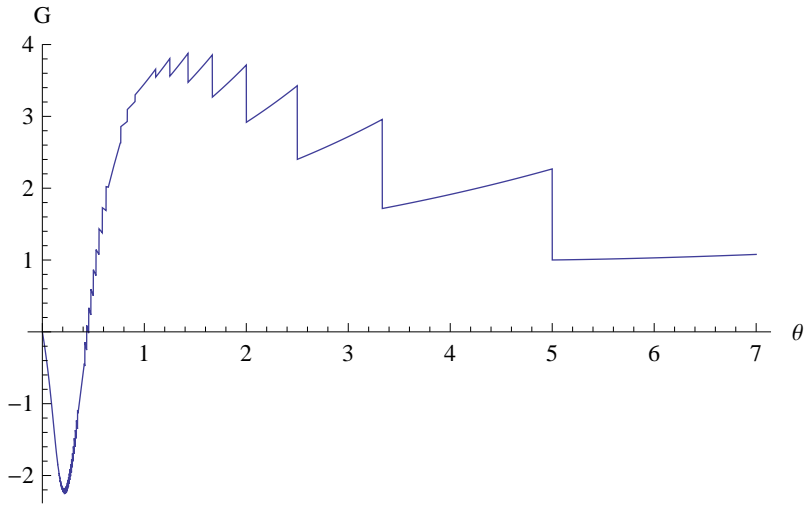}
\caption{$G$ evaluated at $\alpha =1/10$.}
\label{Fig3}
\end{center}
\end{figure}

\begin{figure}[H]
\begin{center}
\includegraphics[scale=1]{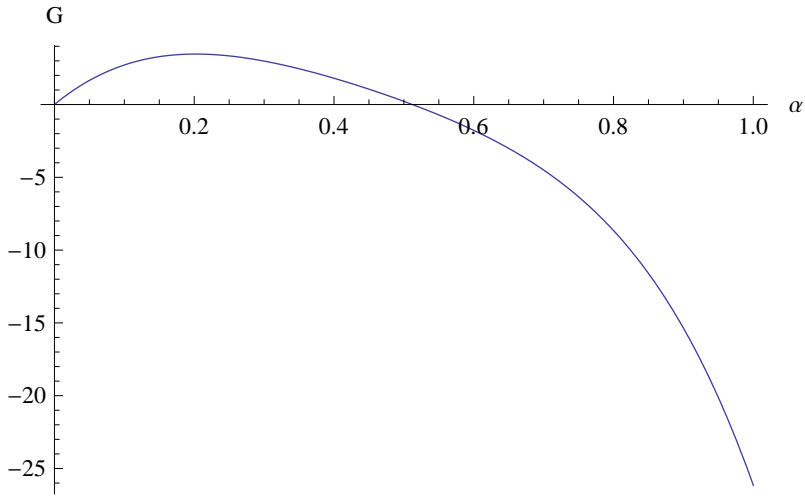}
\caption{$G$ evaluated at $\theta = 3$.}
\label{Fig4}
\end{center}
\end{figure}

We can also describe what happens to the objective when we vary the speed
with which the agents harvests the renewable resource but fix the harvesting
share. Plotting (\ref{Gd}) for $T=10,\sigma =3/20,y_{0}=1,$ we find the
following plot for $G(\theta ,\alpha =1/10)$ in $\theta $
(Figure~\ref{Fig3})
and for
$G(\theta =3,\alpha )$ in $\alpha $
(Figure~\ref{Fig4}).

Note that for varying speeds, the \textquotedblleft
zigzags\textquotedblright \ of the objective result from the fact that for
certain speeds the agent may just manage to complete the final round so that
there is no \textquotedblleft gap\textquotedblright \ in the distribution of
the resource that results from the agent's initial position. This of course
is specific to the assumption of having a uniform distribution of the
resource, i.e. $y_{0}=1$ and relaxed in Ani\c{t}a et.~al.~(2016).

\begin{figure}[H]
\begin{center}
\includegraphics[scale=0.7]{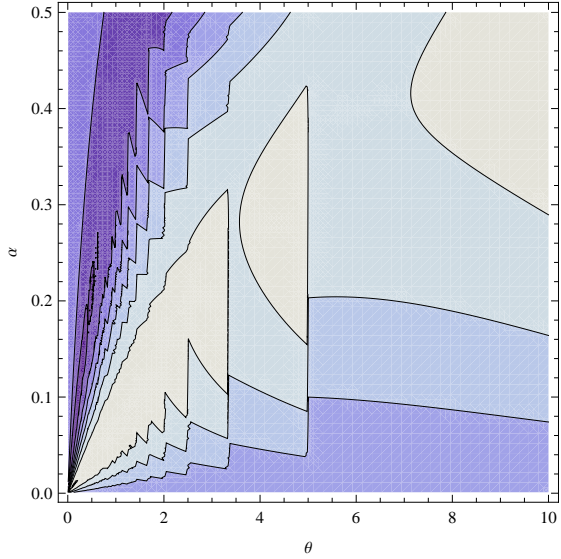}
\caption{Objective function $G$.} 
\label{Fig6}
\end{center}
\end{figure}

Also we present a general contour plot for $G(\theta ,\alpha )$ in
$\alpha, \theta$ (Figure~\ref{Fig6})
where we allow for the harvesting share and the speed
of the agent to vary independently. Note that higher values of the objective
function are characterized by lighter colours. The \textquotedblleft
zigzags\textquotedblright \ carry over into the picture as they are a
property of varying speed only for any chosen value of the harvesting share.

The comparative statics reveals that higher speed implies a larger optimal
exploitation, as the agent will have an increased concern for growth at any
particular location for more frequent returns. Given that the speed is such
that the agent completes exactly one round, the optimal solution coincides
with the static monopoly problem.

\section{Durable good analysis}

In this section we explore the case of a durable good. With the good being
durable, the quantity bought by consumers does not perish and may thus be
consumed later. In this way, the amount of the commodity supplied to (and
sold on) the market decreases demand in later periods - thus intertemporal
demand effects result. Assume that the agent selects speed so that he/she
accomplishes to complete $N$ full rounds of circling and harvesting, i. e. $%
\theta =T/N$. Then, the new \emph{net} revenue function takes the form%
\begin{equation*}
R(\alpha (n))=\alpha (n)\times \left( 1-\sum^{n}\alpha (n)\right)
\end{equation*}%
so that earlier supplies to the market will decrease the marginal return on
later ones. The present value from the $n$th arrival at location $x$ is then
\begin{equation*}
y_{0}\alpha _{n}(1-\sum_{i=1}^{n}\alpha _{i})e^{t_{n}(r-\rho
)}\prod_{i=1}^{n-1}(1-\alpha _{i}).
\end{equation*}%
As we sum over $N$ periods or circling rounds at constant speed $v$ we have
\begin{equation}
\begin{aligned}
\tilde{G}(\alpha)=&\sum_{n=1}^{N}y_{0}\alpha_{n}(1-\sum_{i=1}^{n}\alpha
_{i})e^{t_{n}(r-\rho)}\prod_{i=1}^{n-1}(1-\alpha _{i})\\ =&y_{0}\alpha
_{1}(1-\alpha _{1})e^{t_{1}(r-\rho )}\\ &+y_{0}\alpha _{2}(1-(\alpha
_{1}+\alpha _{2}))e^{t_{2}(r-\rho )}(1-\alpha_{1})\\ &+y_{0}\alpha
_{3}(1-(\alpha _{1}+\alpha _{2}+\alpha _{3}))e^{t_{3}(r-\rho)}(1-\alpha
_{1})(1-\alpha _{2})\\ &+...\\ &+y_{0}\alpha _{N}(1-\sum_{i=1}^{N}\alpha
_{i})e^{t_{N}(r-\rho)}\prod_{i=1}^{N-1}(1-\alpha _{i}). \end{aligned}
\label{T}
\end{equation}%
Clearly the spatial dimension of the problem remains relevant as the agent
returns to any position in future rounds. We can now show that:

\begin{proposition}
Let $K=\{ \alpha =\left( \alpha _{1},\alpha _{2},\dots ,\alpha _{N}\right)
\in \mathbb{R}^{N};\;0\leq \alpha _{j},$ $\forall j\in \{1,2,\dots ,N\},$ $%
\sum \limits_{j=1}^{N}\alpha _{j}\leq 1\}$. Then $\tilde{G}(\alpha )$
attains a global maximum in $\alpha ^{\ast }=\left( \alpha _{1}^{\ast
},\alpha _{2}^{\ast },\dots ,\alpha _{N}^{\ast }\right) ~$on $K$. Also there
are only two situations: I) $\alpha ^{\ast }\in Int(K),$ II) $\alpha
_{1}^{\ast }=\alpha _{2}^{\ast }=\dots =\alpha _{n}^{\ast }=0$ and $\alpha
_{n+1}^{\ast }\neq 0,\dots ,\alpha _{N}^{\ast }\neq 0$.
\end{proposition}

\begin{proof}
See Appendix.
\end{proof}

Examples for $N=2$ and $N=3$ in the appendix show that the maximum of $%
\tilde{G}(\alpha )$ is attained for $\alpha ^{\ast }\in Int(K)$. We further
find that:

\begin{lemma}
With slow growth we get \emph{Cournot type solutions} for $N$ rounds of the
form $\alpha ^{\ast }=\left( \alpha _{1}^{\ast },\alpha _{2}^{\ast },\dots
,\alpha _{N}^{\ast }\right) $, with
\begin{equation*}
\alpha _{j}^{\ast }\approx \frac{1}{N+1}.
\end{equation*}
\end{lemma}

\begin{proof}
See Appendix.
\end{proof}

We therefore generally find a convergence of $\alpha $ in $o(1/N)$. This
result is intuitive as with a durable, spatially heterogeneous renewable
resource commodity, the monopolist plays against itself each round and thus
against its own time-variant copies. With the monopoly result at hand, we
now turn to the case of a non-cooperative game between a finite number of
players.

\section{A durable good game}

Denote a normal form game by $\Gamma =(I,X,u)$ were $i=\left \{
1,2,...,I\right \} $ is the set of symmetric players with (finite)
strategies $X_{i}$ from the set of strategy profiles $A=\displaystyle \times
_{i=1}^{I}A_{i}$. Then $u_{i}:X\rightarrow \mathbb{R}^{+}$ is the payoff
function for player $i$ and $u=(u_{1},...u_{I})$ the payoff function of the
game. The new individual \emph{net} revenue functions now take on the game
form%
\begin{equation*}
R(\alpha (n),i)=\alpha (n,i)\times \left( 1-\sum^{n}\sum^{i}\alpha
(n,i)\right)
\end{equation*}%
For a fixed location and $N$ periods we therefore have individual payoffs
given by:

\begin{equation}
u_{i}(\alpha _{-i})=y_{0}\sum_{n=1}^{N}\alpha
_{i,n}(1-\sum_{j=1}^{I}\sum_{y=1}^{n}\alpha _{j,y})e^{t_{n}(r-\rho
)}\prod_{i=1}^{n-1}(1-\alpha _{i})\text{ }\forall i=\left \{ 1,...,I\right
\} .  \label{Ggame}
\end{equation}

A \emph{Nash equilibrium} of $\Gamma $ is a strategy profile $\alpha ^{\ast
}=(\alpha _{i}^{\ast },\alpha _{-i}^{\ast })$ such that for any player $i$
\begin{equation*}
u_{i}(\alpha ^{\ast })\geq u_{i}(\alpha _{i},\alpha _{-i}^{\ast }), \quad
\forall \alpha _{i}\in A_{i}.
\end{equation*}

\begin{lemma}
 The durable good game with $I$ firms and slow growth has a
symmetric Nash equilibrium%
\begin{equation*}
\alpha _{i}^{\ast }\approx \frac{1}{I(N-1)+2}\quad \forall i=\left \{
1,...,I\right \} .
\end{equation*}
\end{lemma}

\begin{proof}
See Appendix.
\end{proof}

\subsection{An instructive durable good game example}

What happens if growth is not small? We now present an instructive example
for a game with two rounds and an arbitrary number of players.

For two rounds, $N=2,$ and a total number of $I$ players that put $\beta _{j}
$ each in round $j$ in on the market, we have the objective as
\begin{equation*}
\begin{aligned} \left. \tilde{G}(\alpha )\right \vert _{N=2} =\,
&y_{0}\left( \alpha _{1}(1-\alpha _{1}-(I-1)\beta _{1})e^{t_{1}}+\alpha
_{2}\left(1-\left(\alpha _{1}+\alpha _{2}\right.\right.\right.\\
&\left.\left.\left.+(I-1)(\beta _{1}+\beta
_{2})\right)\right)e^{t_{2}}(1-\alpha _{1})\right) \end{aligned}
\end{equation*}%
where again we neglect discounting to avoid clutter. From the foc $\partial
\tilde{G}(\alpha )/\partial \alpha _{1}=0$ we find
\begin{equation*}
\alpha _{1}=\frac{e^{t_{1}}+\beta _{1}e^{t_{1}}-2\alpha _{2}e^{t_{2}}+\alpha
_{2}^{2}e^{t_{2}}-\alpha _{2}\beta _{1}e^{t_{2}}-\alpha _{2}\beta
_{2}e^{t_{2}}-\beta _{1}Ie^{t_{1}}+\alpha _{2}\beta _{1}Ie^{t_{2}}+\alpha
_{2}\beta _{2}Ie^{t_{2}}}{2e^{t_{1}}-2\alpha _{2}e^{t_{2}}}
\end{equation*}%
From the foc $\partial \tilde{G}(\alpha )/\partial \alpha _{2}=0$ we find
\begin{equation*}
\alpha _{2}=\frac{1}{2}\left( 1-\alpha _{1}-\beta _{1}(I-1)-\beta
_{2}(I-1)\right)
\end{equation*}%
for any growth pattern. Solving simultaneously yields optimal shares%
\begin{equation*}
\alpha _{1}=\frac{1}{3e^{t_{2}}}\left( -4e^{t_{1}}+3e^{t_{2}}+\sqrt{%
\begin{array}{c}
4e^{t_{1}}(4e^{t_{1}}-3e^{t_{2}})+ \\
\allowbreak (I-1)^{2}\left( \beta _{1}+\beta _{2}\right) ^{2}e^{2t_{2}}+ \\
\allowbreak 4(I-1)\left( \beta _{1}+4\beta _{2}\right) e^{t_{2}}e^{t_{1}}%
\end{array}%
}-2(I-1)(\beta _{1}+\beta _{2})e^{t_{2}}\right)
\end{equation*}%
and%
\begin{equation*}
\alpha _{2}=\frac{1}{2}\left( 1-\alpha _{1}^{\ast }-(I-1)(\beta _{1}+\beta
_{2})\right) .
\end{equation*}%
Note that for $I=1$ we obtain the results from the durable good monopoly
example.

\begin{figure}
\begin{center}
\includegraphics[scale=0.5]{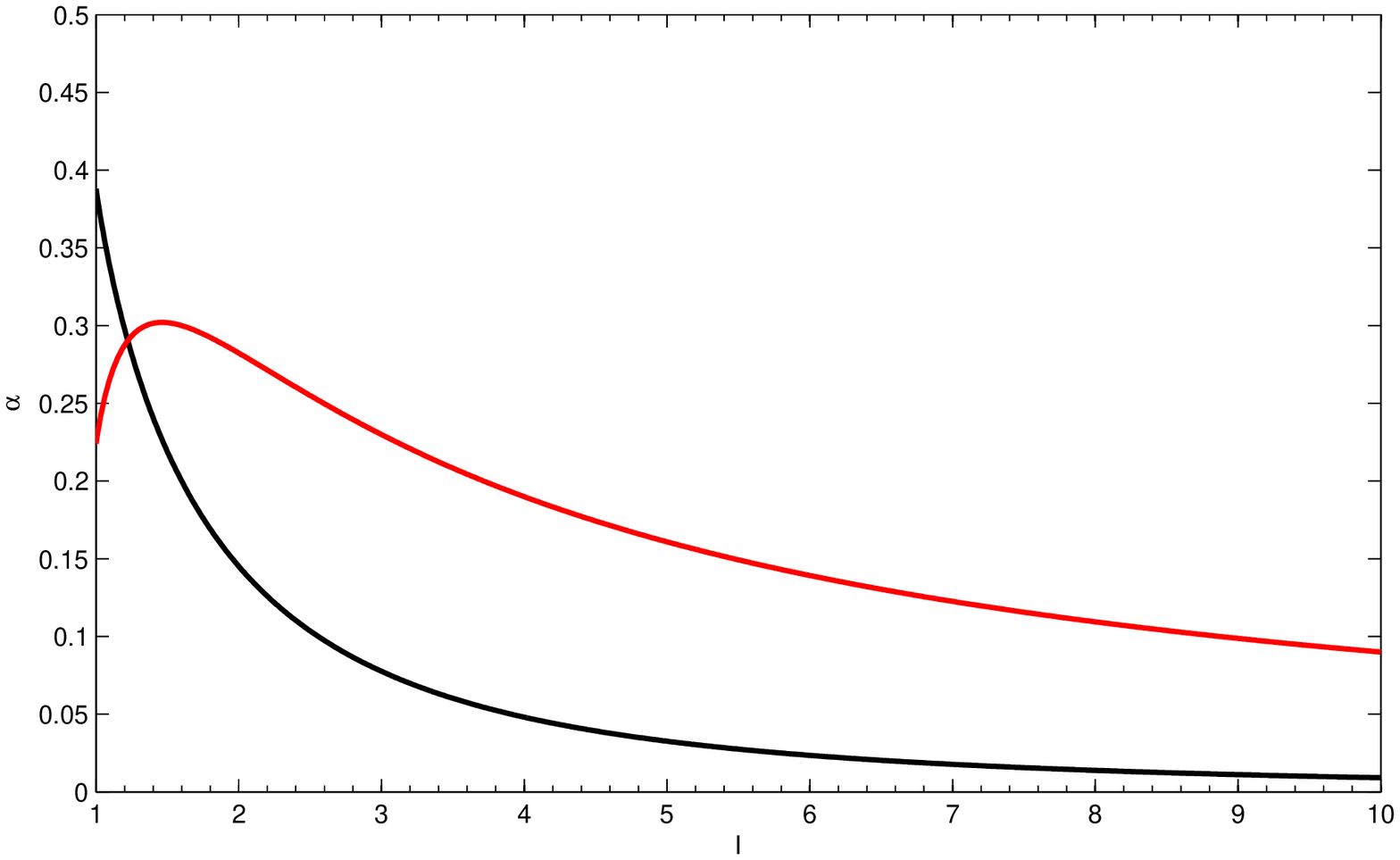}
\caption{Game figure.}
\label{Fig7}
\end{center}
\end{figure}

Assume\emph{\ }growth again satisfies $\frac{3}{4}e^{t_{2}}<e^{t_{1}}<e^{t2}$
as in the durable good analysis above. In particular let $e^{t_{1}}=1$ and $%
e^{t_{2}}=11/9.$ To solve for symmetric\ equilibrium over rounds, the
optimal share $\alpha _{1}$ of a player needs to be the same as the optimal
share as that of all the other players $\beta .$ This is the case in both
rounds. Hence with $\beta _{1}=\alpha _{1},\beta _{2}=\alpha _{2},$ solving
the resulting system simultaneously again, we find the equilibrium share in
round $1$ with $I$ players as:%
\begin{equation*}
\alpha _{1}^{\ast }=\frac{1}{I}-(\frac{1}{I}\frac{I+1}{44I+22}\left(
9I^{2}+7I+20-\sqrt{81I^{4}+126I^{3}+409I^{2}-512I+4}\right)
\end{equation*}%
and the equilibrium share in round $2$ with $I$ players as:%
\begin{equation*}
\alpha _{2}^{\ast }=\frac{1}{44I+22}\left( 9I^{2}+7I+20-\sqrt{%
81I^{4}+126I^{3}+409I^{2}-512I+4}\right)
\end{equation*}%
Both can be graphically represented, with round $1$ in red, see
Figure~\ref{Fig7}.

\begin{figure}[H]
\begin{center}
\includegraphics[scale=0.5]{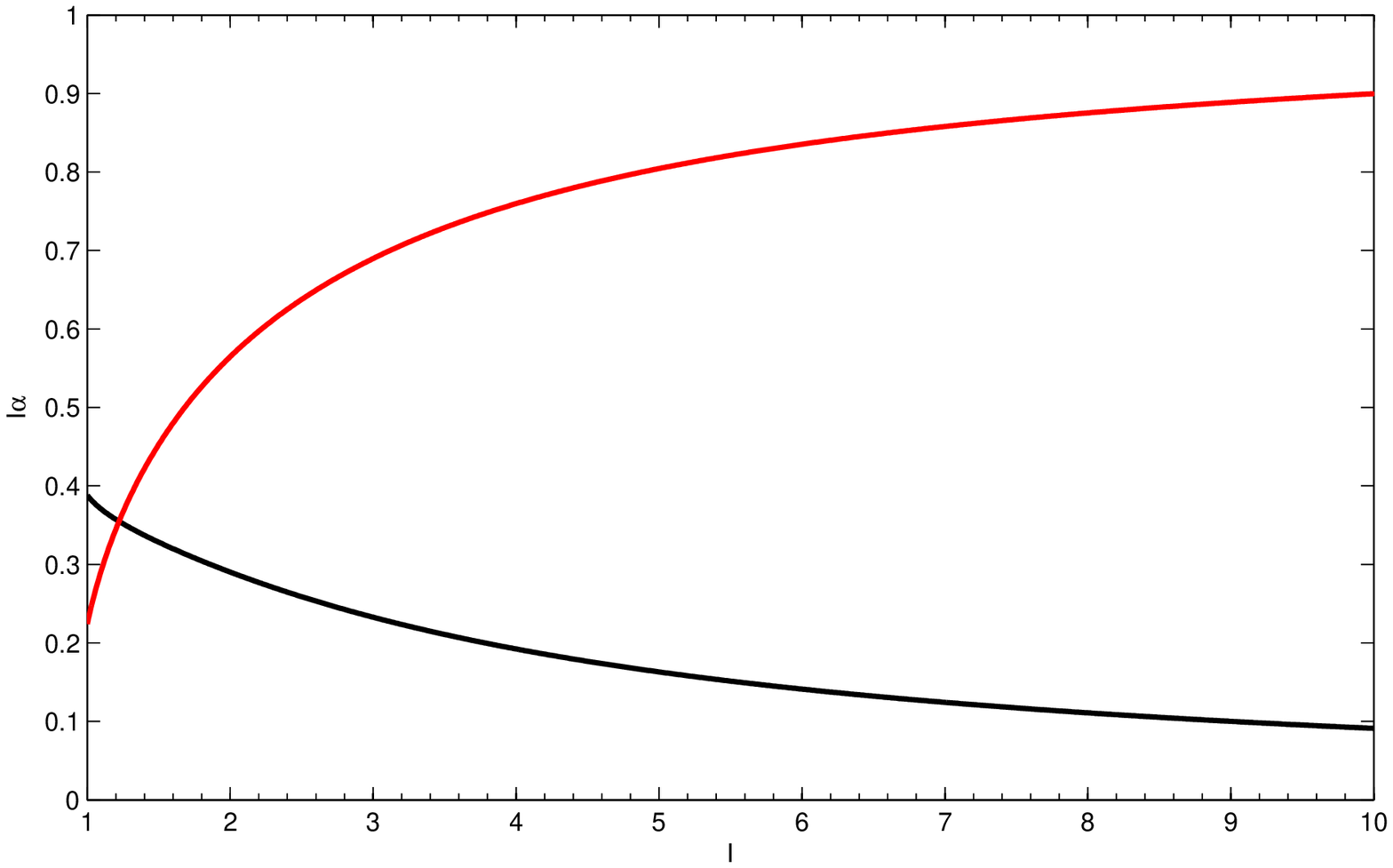}
\caption{Game figure}
\label{Fig8}
\end{center}
\end{figure}

Confirming our earlier results, for $I=1$ (\emph{monopoly}) we find $\alpha
_{1}^{\ast }\approx 0.224\,$ and $\alpha _{2}^{\ast }\approx 0.388\,,$
leaving more than half of the resource unharvested. For \emph{duopoly} we
find for each player equilibrium shares of $\alpha _{1}^{\ast }\approx
0.282\, \,$\ and $\alpha _{2}^{\ast }\approx \allowbreak 0.145.$ It can be
shown straightforwardly that payoffs (profits) $\left. \tilde{G}(\alpha
)\right \vert _{N=2}$ are monotonically decreasing in $I$, as is the final
resource stock. The limits for the total harvested shares each round ($1$ in
red again) are given in Figure~\ref{Fig8}.

Note that the resource is \emph{fully depleted asymptotically already in
round }$1$ (even without discounting) so that competition between more
players is detrimental for the resource both in a quantity and a time
dimension. Hence with many players we find an equilibrium result in this
game that is very different from the cases I and II in monopoly.

Generally we observe that for the equilibrium strategies convergence
satisfies
\begin{equation*}
\alpha _{i}^{\ast }=o\left( 1/(NI)\right) .
\end{equation*}%
The symmetric equilibrium limit price (with many players $I$ and/or many
rounds $N$) for a given demand is thus:%
\begin{equation*}
\begin{aligned} \lim_{N\rightarrow \infty }p^{\ast } &=\lim_{\substack{
N\rightarrow \infty \\ \text{or }I\rightarrow \infty }}\left(
1-\sum_{y=1}^{I}\sum_{j=0}^{N}\alpha _{j,y}\left(t+\theta j\right)\right) \\
&=\lim_{\substack{ N\rightarrow \infty \\ \text{or }I\rightarrow \infty
}}\left( 1-NIo(1/NI)\right) = 0, \end{aligned}
\end{equation*}%
which is the normalized level of marginal costs in this model so that we
have \emph{Cournot convergence} to the Walrasian price (perfect competition).

\section{The limit principle}

Green's limit principle is derived in a repeated game setting. With $I$
symmetric players (firms)\ and $N$ \emph{discrete rounds} with fixed
locations we found symmetric equilibrium shares of
\begin{equation*}
\alpha _{i}^{\ast }=\frac{1}{I(N-1)+2}
\end{equation*}%
for each round of harvesting of the durable good with $N,I>1$.

A general problem in \emph{continuous time} games is that looking at
deviations from collusive behaviour is generically trivial, as they are
immediately detected and punished. Here we have a natural way of
\textquotedblleft discretizing\textquotedblright \ the continuum by assuming
that individual deviation detection takes a full round so that deviation
profits for firm $i$ from the shared monopoly outcome are%
\begin{equation*}
\pi _{i}^{d}=\max_{\alpha _{i}\in \mathcal{A}}\int_{0}^{t_{1}}\left(
e^{-\rho t}\alpha _{i}(t)\left( 1-\alpha _{i}(t)-\sum_{j\neq i}\alpha
_{j}^{\ast }(t)\right) f^{\alpha }(t-)-C(\alpha (t))\right) dt.
\end{equation*}%
Under discounting, deviation will take place in the first round or never%
\footnote{%
This changes if we take the growth process into account.}. The deviation
profit satisfies%
\begin{equation*}
\pi _{i}^{d}>\int_{0}^{t_{1}}\left( e^{-\rho t}\alpha _{i}^{\ast }(t)\left(
1-\alpha _{i}^{\ast }(t)-\sum_{j\neq i}\alpha _{j}^{\ast }(t)\right)
f^{\alpha }(t-)-C(\alpha (t))\right) dt\equiv \pi _{i}^{\ast }
\end{equation*}%
i.e. exceeds the equilibrium payoff during the punishment rounds.

As in Green (1980) we assume that in the replica economy indexed by $r$,
the individual demands are equally scaled up by $I$. Then the individual
firm monopoly level, i.e. the \emph{collusive outcome} with equally shared
production\ is
\begin{equation*}
\alpha _{i}^{M,r}=\frac{1}{N+1}\quad \forall i=\left \{ 1,...,I\right \}
\end{equation*}%
and deviation may be punished (ad infinitum)\ by playing the scaled
symmetric Cournot equilibrium
\begin{equation*}
\alpha _{i}^{\ast ,r}=I\alpha _{i}^{\ast }=\frac{I}{I(N-1)+2}.
\end{equation*}%
Note that $\alpha _{i}^{\ast ,r}>\alpha _{i}^{M,r}$ as $I>1.$

\begin{lemma}
 In the replica economy, there exists a discount rate $\rho >0$
small enough, such that the collusive outcome is an equilibrium in the
continuous time game.
\end{lemma}

\begin{proof}
See Appendix.
\end{proof}

Now assume that individual shares are subject to some\emph{\ idiosyncractic
noise} term$~\tilde{\alpha}_{i}=\alpha _{i}+\tilde{\varepsilon}_{i}$ where $%
\tilde{\varepsilon}_{i}\sim N(0,Var(\varepsilon ))$ with $0<Var(\varepsilon
)<\infty $ and that the cartel can observe only some aggregate statistic of
play, e.g. \emph{the price in the replica economy}, where aggregate demand $%
I(1-p)$ equals aggregate supply $\sum_{i=1}^{I}\alpha _{i}$. This random
variable that satisfies%
\begin{equation*}
\tilde{p}^{r}=1-\frac{1}{I}\sum_{i=1}^{I}\left( \alpha _{i}^{r}+\tilde{%
\varepsilon}_{i}\right) .
\end{equation*}

We then find that:

\begin{proposition}
With idiosyncratic noise $\tilde{\varepsilon}_{i}\sim N(0,Var(\varepsilon
)), $ individual deviations become undetectable and we get the limit
principle to hold in the durable good game.
\end{proposition}

\begin{proof}
See Appendix.
\end{proof}

As observed in Mas-Colell (1988, p.30) this finding is somewhat paradoxical
in the theory of perfect competition as it is \emph{not} perfect information
but noise that helps perfect competition to come about. In the \emph{%
non-durable commodity} case, as has been shown in section 2 and by
Ani\c{t}a et.~al.~(2016).
the monopoly solution often implies letting the resource grow
unimpaired and harvest the profit maximizing quantity only in the last
round. This however is not implementable with more than one player, as
deviation in this last round cannot be punished. Hence the optimal collusive
outcome cannot be sustained even without noise and the limit principle holds.

\section{Conclusion}

In this paper we have shown that endogenizing prices \emph{prevents the
extinction of the renewable resource} compared to the 2014 model of
Behringer and Upmann. Letting prices fluctuate therefore presents an
alternative policy to forcing the agent to go multiple rounds or to move
with a minimum speed in order to make him/her take into account the future
more seriously.

Also we have shown that in a harvesting game, competition will have a
critical temporal dimension in addition to the negative effects on the stock
of the renewable resource in that the resource is depleted earlier. Optimal
shares of the harvested renewable resource in this fully dynamic spatial
model inherit the convergence properties of the static Cournot model. Hence
the Walrasian properties, implying perfect competition with many firms but
also the dynamic results of Green (1980) for a stationary repeated game
setting are reestablished. Green's limit principle is shown to be robust to
the investigation of competition in durable and non-durable renewable
resources in our non-stationary setting.

\section{Appendix}

\paragraph{Proof of Lemma 1.}

Rearranging (\ref{0}) we get
\begin{equation*}
\begin{aligned} 0 &\geq -\sum_{l=0}^{k}\int_{0}^{\theta}e^{-\rho (t+\theta
l)}\alpha ^{\ast }\left(t+\theta l\right)\left( 1-\alpha ^{\ast
}\left(t+\theta l\right)\right) f_{l}^{\alpha ^{\ast }}(t-)dt \\
&+\sum_{l=0}^{k}\int_{0}^{\theta}e^{-\rho (t+\theta l)}(\alpha ^{\ast
}+\varepsilon w)\left(t+\theta l\right)( 1-\alpha ^{\ast }-\varepsilon
w)\left(t+\theta l\right) f_{l}^{\alpha ^{\ast }+\varepsilon w}(t-)dt
\end{aligned}
\end{equation*}%
or
\begin{equation*}
\begin{aligned} 0&\geq \sum_{l=0}^{k}\int_{0}^{\theta}e^{-\rho (t+\theta
l)}\left[ -\alpha ^{\ast }\left(t+\theta l\right)\left( 1-\alpha ^{\ast
}\left(t+\theta l\right)\right) f_{l}^{\alpha ^{\ast }}(t-)\right.\\
&+\left.(\alpha ^{\ast }+\varepsilon w)\left(t+\theta l\right)( 1-\alpha
^{\ast }-\varepsilon w)\left(t+\theta l\right) f_{l}^{\alpha ^{\ast
}+\varepsilon w}(t-) \right]dt. \end{aligned}
\end{equation*}%
Expanding, transforming, and dividing by $\varepsilon >0$ yields
\begin{equation*}
\begin{aligned} 0&\geq \sum_{l=0}^{k}\int_{0}^{\theta}e^{-\rho (t+\theta
l)}\left[ w\left(t+\theta l \right)\left( 1-2\alpha ^{\ast }\left(t+\theta l
\right)\right. \right.\\ &\left.\left.-\varepsilon w\left(t+\theta l
\right)\right) f_{l}^{\alpha ^{\ast }}(t-)+\left(\alpha ^{\ast
}\left(t+\theta l \right)\right.\right.\\ &\left.\left.+\varepsilon
w\left(t+\theta l \right)\right)\left( 1-\alpha ^{\ast }\left(t+\theta l
\right)-\varepsilon w\left(t+\theta l \right)\right)\frac{\left(
f_{l}^{\alpha^{\ast }+\varepsilon w}-f_{l}^{\alpha ^{\ast
}}\right)(t-)}{\varepsilon } \right]dt. \end{aligned}
\end{equation*}%
Taking $\varepsilon \rightarrow 0$ we find%
\begin{equation*}
\begin{aligned} 0&\geq \sum_{l=0}^{k}\int_{0}^{\theta}e^{-\rho (t+\theta
l)}\left[w\left(t+\theta l\right)\left( 1-2\alpha ^{\ast }\left(t+\theta
l\right)\right) f_{l}^{\alpha ^{\ast }}(t-)\right.\\ &\left.+\alpha ^{\ast
}\left(t+\theta l\right)\left( 1-\alpha ^{\ast }\left(t+\theta
l\right)\right) z_{l}(t)\right] dt. \end{aligned}
\end{equation*}%
Noting that
\begin{equation*}
z_{l}=\lim_{\varepsilon \rightarrow 0}\frac{f_{l}^{\alpha ^{\ast
}+\varepsilon w}-f_{l}^{\alpha ^{\ast }}}{\varepsilon }\quad \in \quad
L^{\infty }(0,T)
\end{equation*}%
and using (\ref{den}) we get the conclusion.
\qed

\paragraph{Proof of Proposition 2.}

We multiply the first equation in (\ref{5}) by $z_{l}(t)$, integrate on $%
\left[ 0,\theta \right) $ and add up over $l$ to $k-1$. We get that
\begin{equation*}
\begin{aligned} &\sum_{l=0}^{k-1}\int_{0}^{\theta}p_{l}(t)z_{l}(t)dt
=\sum_{l=0}^{k-1}\int_{0}^{\theta}\left[ e^{r\theta}\left( 1-\alpha ^{\ast
}\left( t+\theta l\right) \right) p_{l+1}(t)z_{l}(t)\right.\\
&\left.+e^{-\rho (t+\theta l)}\alpha ^{\ast }\left( t+\theta l\right) \left(
1-\alpha ^{\ast }\left( t+\theta l\right) \right) z_{l}(t)\right] dt.
\end{aligned}
\end{equation*}%
Now replace from (\ref{4.1}) that
\begin{equation*}
e^{r\theta }\left( 1-\alpha ^{\ast }\left( t+\theta l\right) \right)
z_{l}(t)=z_{l+1}(t)+e^{r\theta }w\left( t+\theta l\right) f_{l}^{\alpha
^{\ast }}(t-),
\end{equation*}%
so that
\begin{equation*}
\begin{aligned}
&\sum_{l=0}^{k-1}\int_{0}^{\theta}p_{l}(t)z_{l}(t)dt=\sum_{l=0}^{k-1}%
\int_{0}^{\theta}p_{l+1}(t)\left[z_{l+1}(t) +e^{r\theta}w\left(t+\theta
l\right)f_l^{\alpha^{\ast}}(t-)\right] dt\\
&+\sum_{l=0}^{k-1}\int_{0}^{\theta}e^{-\rho (t+\theta l)}\alpha ^{\ast
}\left( t+\theta l\right) \left( 1-\alpha ^{\ast }\left( t+\theta l\right)
\right) z_{l}(t) dt. \end{aligned}
\end{equation*}%
Since $z_{0}(t)=0$ and $p_{k}(t)$ satisfies the second equation in (\ref{5}%
), we may conclude that
\begin{equation*}
\begin{aligned} &\sum_{l=0}^{k-1}\int_{0}^{\theta}e^{r\theta}w\left(t+\theta
l\right)f_{l}^{\alpha ^{\ast }}(t-)p_{l+1}(t))dt \\ &+\int_{0}^{T-k\theta}
e^{-\rho (t+\theta k)}\alpha ^{\ast }\left(t+\theta k\right)\left( 1-\alpha
^{\ast }\left(t+\theta k\right)\right) z_{k}(t)dt \\
&+\sum_{l=0}^{k-1}\int_{0}^{\theta} e^{-\rho (t+\theta l)}\alpha ^{\ast
}\left(t+\theta l\right)\left( 1-\alpha ^{\ast }\left(t+\theta
l\right)\right) z_{l}(t) dt=0 \end{aligned}
\end{equation*}%
and
\begin{equation}
\begin{aligned} &\sum_{l=0}^{k}\int_{0}^{\theta} e^{-\rho (t+\theta
l)}\alpha ^{\ast }\left(t+\theta l\right)\left( 1-\alpha ^{\ast
}\left(t+\theta l\right)\right) z_{l}(t) dt\\
&=-\sum_{l=0}^{k-1}\int_{0}^{\theta} e^{r\theta}w\left(t+\theta
l\right)f_{l}^{\alpha ^{\ast }}(t-)p_{l+1}(t))dt \end{aligned}  \label{s}
\end{equation}%
Using (\ref{3}) we obtain
\begin{equation}
\begin{aligned} 0 \geq &\sum_{l=0}^{k-1}\int_{0}^{\theta}\left[ w\left(
t+\theta l\right) f_{l}^{\alpha ^{\ast }}(t-) e^{-\rho (t+\theta l)}\left(
1-2\alpha ^{\ast }\left( t+\theta l\right) \right)\right. \\
&\left.+e^{-\rho (t+\theta l)}\alpha ^{\ast }\left( t+\theta l\right) \left(
1-\alpha ^{\ast }\left( t+\theta l\right) \right) z_{l}(t)\right] dt \\
&+\int_{0}^{\theta}e^{-\rho (t+\theta k)}w\left( t+\theta k\right) \left(
1-2\alpha ^{\ast }\left( t+\theta k\right) \right) f_{k}^{\alpha ^{\ast
}}(t-)dt. \end{aligned}
\end{equation}%
From (\ref{s}) we get now that
\begin{equation}
\begin{aligned} 0 \geq &\sum_{l=0}^{k-1}\int_{0}^{\theta}w\left(t+\theta
l\right)\left[ e^{-\rho (t+\theta l)}\left( 1-2\alpha ^{\ast }\left(t+\theta
l\right)\right)\right.\\ &\qquad \qquad \left.-e^{r\theta}p_{l+1}(t)\right]
f_{l}^{\alpha ^{\ast }}(t-)dt \\ &+\int_{0}^{\theta}e^{-\rho (t+\theta
k)}w\left(t+\theta k\right)\left( 1-2\alpha ^{\ast }\left(t+\theta
k\right)\right) f_{k}^{\alpha ^{\ast }}(t-)dt, \end{aligned}
\end{equation}%
$\forall \,w$ (which is arbitrary) such that $0\leq \alpha ^{\ast
}+\varepsilon w\leq 1,\,a.e$.
\qed

\paragraph{Proof of Proposition 3.}

\emph{One} round of cycling yields
\begin{equation*}
\begin{aligned} E(1)=&\alpha (1-\alpha )y_{0}\int_{0}^{2\pi }e^{(r-\rho
)t_{1}(x)}dx=\alpha(1-\alpha ) y_{0}\int_{0}^{2\pi }e^{(r-\rho )\frac{\theta
x}{2\pi }}dx\\ =&\alpha(1-\alpha ) y_{0}\frac{2\pi }{(r-\rho )\theta }\left(
e^{(r-\rho )\theta }-1\right), \end{aligned}
\end{equation*}%
so total discounted supply in the $n$th period is%
\begin{equation*}
\begin{aligned} E(n) =&\alpha \left( 1-\sum_{i=1}^{n}\alpha _{i}\right)
(1-\alpha )^{n-1}y_{0}\int_{0}^{2\pi }e^{-\rho t_{n}(x)}e^{\left(
(n-1)r\theta +r\frac{\theta x}{2\pi }\right) }dx \\ =&\alpha \left(
1-n\alpha \right) (1-\alpha )^{n-1}\frac{2\pi y_{0}}{(r-\rho )\theta }\left(
e^{(r-\rho )\theta }-1\right) e^{(n-1)(r-\rho )\theta }. \end{aligned}
\end{equation*}%
Summing over all periods (defining the net growth rate as $\sigma \equiv
r-\rho )$ this yields%
\begin{eqnarray*}
\sum_{i=1}^{n} E(i) &=&
\sum_{i=1}^{n} \alpha (1-i\alpha )(1-\alpha )^{i-1}%
\frac{2\pi y_{0}}{\sigma \theta }(e^{\sigma \theta }-1)e^{(i-1)\sigma \theta
} \\
&=&-\frac{2\pi y_{0}}{\sigma \theta }\frac{\alpha }{(\alpha +e^{-\sigma
\theta }-1)^{2}}\times \\
&&\left(
\begin{array}{c}
(\alpha -1)(2e^{-\sigma \theta }-e^{-2\sigma \theta }-1) \\
+\left(
\begin{array}{c}
(\alpha -1+n\alpha )e^{\sigma \theta (n+1)}+ \\
(n\alpha ^{2}-2\alpha +2-2n\alpha )e^{\sigma \theta (n+2)}+ \\
(\alpha -1-n\alpha ^{2}+n\alpha )e^{\sigma \theta (n+3)}%
\end{array}%
\right) e^{-3\sigma \theta }(1-\alpha )^{n}%
\end{array}%
\right)
\end{eqnarray*}%
if $(1-\alpha )e^{\sigma \theta }-1\neq 0.$

For the possibly incomplete final round we have
\begin{equation*}
E(n)=\alpha \left( 1-n\alpha \right) (1-\alpha )^{n-1}\frac{2\pi y_{0}}{%
\sigma \theta }\left( e^{\sigma \theta }-1\right) e^{(n-1)\sigma \theta }
\end{equation*}%
so that
\begin{equation*}
E(n,x)=\alpha \left( 1-n\alpha \right) \alpha ^{n-1}\frac{2\pi y_{0}}{\sigma
\theta }\left( e^{\sigma \frac{\theta x}{2\pi }}-1\right) e^{(n-1)\sigma
\theta }
\end{equation*}%
and%
\begin{equation*}
E(n+1,s(T))=\alpha \left( 1-(n+1)\alpha \right) \alpha ^{n}\frac{2\pi y_{0}}{%
\sigma \theta }\left( e^{\sigma \frac{\theta s(T)}{2\pi }}-1\right)
e^{n\sigma \theta }
\end{equation*}%
So if we add up to $N$ we get%
\begin{equation*}
G(\theta ,\alpha )=\sum_{i=1}^{N}E(i)+E(N+1,s(T))
\end{equation*}%
as given above.
\qed

\paragraph{Proof of Proposition 4.}

As
\[
K=\left \{ \alpha =\left( \alpha _{1},\alpha _{2},\dots ,\alpha
_{N}\right) \in \mathbb{R}^{N};\;0\leq \alpha _{j},\; \forall j\in
\{1,2,\dots ,N\},\; \sum \limits_{j=1}^{N}\alpha _{j}\leq 1\right \},
\]
we have $\tilde{G}:K\rightarrow \mathbb{R}$. Since $\tilde{G}$ is a continuous
function and $K$ is compact, by Weierstass Theorem we conclude that $\tilde{G%
}$ attains its global maximum on $K$ in $\alpha ^{\ast }=\left( \alpha
_{1}^{\ast },\alpha _{2}^{\ast },\dots ,\alpha _{N}^{\ast }\right) $. We
prove that we have only two situations:

\begin{enumerate}
\item[I.] $\alpha^{*} \in Int(K)$;

\item[II.] $\alpha _{1}^{\ast }=\alpha _{2}^{\ast }=\dots =\alpha _{n}^{\ast
}=0$ and $\alpha _{n+1}^{\ast }\neq 0,\dots ,\alpha _{N}^{\ast }\neq 0$.
We denote by
\[ K_{n} := \left \{ \left( \alpha _{n+1},\dots ,\alpha _{N}\right)
\in \mathbb{R}^{N-n};\;0\leq \alpha _{j},\; \forall j\in \{n+1,\dots ,N\},\;
\sum \limits_{j=n+1}^{N}\alpha _{j}\leq 1\right \}.
\]
Then $\left( \alpha
_{n+1}^{\ast },\alpha _{n+2}^{\ast },\dots ,\alpha _{N}^{\ast }\right) \in
Int(K_{n})$.
\end{enumerate}

We argue by reductio ad absurdum: Assume that $\alpha ^{\ast }\in \partial K$
(where $\partial K$ is the boundary of $K$).

If $\sum \limits_{j=1}^N \alpha_j^{*}= 1$, then consider the smallest $n$
such that $\sum \limits_{j=1}^n \alpha_j^{*}= 1$. It follows that $%
\alpha_{n+1}^ {*}= \alpha_{n+2}^{*}= \dots =0$ and that
\begin{equation*}
\tilde{G}(\alpha^{*})=y_0\sum_{j=1}^{n}\alpha_{j}^{*}(1-\sum_{k=1}^{j}\alpha
_{k}^{*})e^{t_{j}(r-\rho)}\prod_{i=1}^{j-1}(1-\alpha _{i}^{*}).
\end{equation*}
Since $\sum \limits_{j=1}^n \alpha_j^{*}= 1$, if we take $\tilde{\alpha}%
=\left(\alpha_1^{*}, \alpha_2^{*}, \dots, \alpha_{n-1}^{*}, \frac{%
\alpha_{n}^{*}}{2},0,0, \dots, 0\right) \in K$, then
\begin{equation*}
\tilde{G}(\tilde{\alpha})>\tilde{G}(\alpha^{*})
\end{equation*}
and this is a contradiction.

If there is a $n\in \{1,2,\dots ,N\}$ such that $\alpha _{n}^{\ast }=0$,
then consider the smallest such $n$. For $n\geq 2$ and a $\tilde{\alpha}\in
K $ that differs from $\alpha ^{\ast }$ by two components: $\tilde{\alpha}%
_{n-1}=\alpha _{n}^{\ast }=0$ and $\tilde{\alpha}_{n}=\alpha _{n-1}^{\ast
}>0 $, we have
\begin{equation*}
\tilde{G}(\tilde{\alpha})>\tilde{G}(\alpha ^{\ast })
\end{equation*}%
which is a contradiction.

If $\alpha _{1}^{\ast }=\alpha _{2}^{\ast }=\dots =\alpha _{n}^{\ast }=0$
and $\alpha _{n+1}^{\ast }\neq 0,\dots ,\alpha _{N}^{\ast }\neq 0$ we obtain
that $\left( \alpha _{n+1}^{\ast },\alpha _{n+2}^{\ast },\dots ,\alpha
_{N}^{\ast }\right) \in Int(K_{n})$.

We may conclude now that $\alpha ^{\ast }$ is in one of two situations and
so $\alpha ^{\ast }$ is one of the steady states for $\tilde{G}(\alpha )$.
Corresponding to the two cases, $\alpha ^{\ast }$ is the solution of one of
the two systems

\begin{enumerate}
\item[I.] $\displaystyle \frac{\partial \tilde{G}(\alpha )}{\partial \alpha
_{j}} = 0, \quad j \in \{1, 2, \dots, N\}$;

\item[II.] $\left \{
\begin{aligned}
&\alpha_1= \dots =\alpha_n=0, \quad \alpha_{n+1}\neq 0, \dots, \alpha_N\neq 0,\\
&\frac{\partial \tilde{G}(\alpha )}{\partial \alpha _{j}} = 0, \quad j \in \{n+1, \dots, N\}.\\
\end{aligned}\right. $
\end{enumerate}
\qed

\paragraph{Proof of Lemma 5.}

By calculus, we get from (\ref{T}) for the particular case $n=1$ and $%
y_{0}=1 $ that
\begin{equation*}
\begin{aligned} \frac{\partial \tilde{G}(\alpha )}{\partial \alpha
_{1}}&=\left((1-\alpha_{1})+\alpha _{1}(-1)\right) e^{t_{1}(r-\rho )}\\
&+\alpha _{2}\left( (-1)(1-\alpha _{1})+(1-(\alpha
_{1}+\alpha_{2}))(-1)\right) e^{t_{2}(r-\rho )}\\ &+\alpha _{3}\left(
(-1)(1-\alpha _{1})(1-\alpha _{2})+(1-(\alpha_{1}+\alpha _{2}+\alpha
_{3}))(-1)(1-\alpha _{2})\right) e^{t_{3}(r-\rho )}.\\ &+\alpha _{N}\left(
(-1)\prod_{i=1}^{N-1}(1-\alpha_{i})+(1-\sum_{i=1}^{N}\alpha
_{i})(-1)\prod_{i=2}^{N-1}(1-\alpha_{i})\right) e^{t_{N}(r-\rho )} \\
&=(1-\alpha _{1})e^{t_{1}(r-\rho)}+\sum_{k\geq 1}\alpha
_{k}(-1)e^{t_{k}(r-\rho)}\prod_{i=1}^{k-1}(1-\alpha _{i})\\ &+\sum_{k\geq
2}\alpha
_{k}\left(1-\sum_{i=1}^{k}\alpha_i\right)e^{t_{k}(r-\rho)}(-1)\prod_{%
\substack{i=2}}^{k-1}(1-\alpha _{i}).\\ \end{aligned}
\end{equation*}%
For the general case, when $n=1,2,\dots ,N$, we observe that
\begin{equation*}
\begin{aligned} \frac{\partial \tilde{G}(\alpha )}{\partial \alpha _{n}}
&=(1-\sum_{i=1}^{n}\alpha _{i})e^{t_{n}(r-\rho)}\prod_{i=1}^{n-1}(1-\alpha
_{i})\\ &+\sum_{k\geq n}\alpha
_{k}(-1)e^{t_{k}(r-\rho)}\prod_{i=1}^{k-1}(1-\alpha _{i})\\ &+\sum_{k\geq
n+1}\alpha
_{k}\left(1-\sum_{i=1}^{k}\alpha_i\right)e^{t_{k}(r-\rho)}(-1)\prod_{%
\substack{i=2\\i \neq n}}^{k-1}(1-\alpha _{i}). \end{aligned}
\end{equation*}%
For $n=N$ we have the necessary optimality condition for the last round as:%
\begin{equation*}
\frac{\partial \tilde{G}(\alpha )}{\partial \alpha _{N}}=y_{0}(1-%
\sum_{i=1}^{N-1}\alpha _{i}-2\alpha _{N})e^{t_{N}(r-\rho
)}\prod_{i=1}^{N-1}(1-\alpha _{i})=0
\end{equation*}%
Thus for slow growth we have%
\begin{equation*}
\alpha _{i}\approx \alpha _{N}\approx \frac{1}{N+1}
\end{equation*}%
\qed

\paragraph{Proof of Lemma 6.}

The game now has the payoffs as in (\ref{Ggame}). Taking the derivative
w.r.t to the last round, assuming the other firms are symmetric we find a
necessary condition for optimal shares to satisfy:
\begin{eqnarray*}
\frac{\partial u_{i}(\alpha _{-i})}{\partial \alpha _{N}} &=&y_{0}\alpha
_{N}(-1)e^{t_{N}(r-\rho )}\prod_{i=1}^{N-1}(1-\alpha _{i})+ \\
&&y_{0}(1-\left( \sum_{i=1}^{N-1}\alpha _{i}+\alpha
_{N}+(I-1)\sum_{y=1}^{n}\alpha _{j,y}\right) e^{t_{N}(r-\rho
)}\prod_{i=1}^{N-1}(1-\alpha _{i})=0
\end{eqnarray*}%
or%
\begin{equation*}
1-\left( \sum_{i=1}^{N-1}\alpha _{i}+2\alpha
_{N}+(I-1)\sum_{y=1}^{N-1}\alpha _{j,y}\right) =0
\end{equation*}%
or with \emph{symmetric} shares%
\begin{equation*}
\alpha _{i}=\frac{1-((I-1)\sum_{y=1}^{N-1}\alpha _{j,y})}{N+1}
\end{equation*}%
or as then also $\alpha _{i}\approx \alpha _{j,y}$ for a low discounting
rate and we have
\begin{equation*}
\alpha _{i}^{\ast }\approx \frac{1}{I(N-1)+2}
\end{equation*}%
$\forall i=\left \{ 1,...,I\right \} $
\qed

\paragraph{Examples for durable good analysis.}

1. For $N=2$ we then have
\begin{equation*}
\tilde{G}(\alpha )=y_{0}\left( \alpha _{1}(1-\alpha _{1})e^{t_{1}}+\alpha
_{2}(1-(\alpha _{1}+\alpha _{2}))e^{t_{2}}(1-\alpha _{1})\right).
\end{equation*}
We set $y_{0}=1$. $\alpha^{*}$ is the solution of one of the two system

\begin{enumerate}
\item[I.] $\left \{
\begin{aligned}
&\frac{\partial \tilde{G}(\alpha )}{\partial \alpha _{1}} = 0,\\
&\frac{\partial \tilde{G}(\alpha )}{\partial \alpha _{2}} = 0.\\
\end{aligned}
\right. $

\item[II.] $\left \{
\begin{aligned}
&\alpha_1=0, \quad \alpha_{2}\neq 0,\\
&\frac{\partial \tilde{G}(\alpha )}{\partial \alpha _{2}} = 0.\\
\end{aligned}
\right. $
\end{enumerate}

I. From
\begin{equation*}
\frac{\partial \tilde{G}(\alpha )}{\partial \alpha _{1}}=e^{t_{1}}-2\alpha
_{1}e^{t_{1}}-2\alpha _{2}e^{t_{2}}+\alpha _{2}^{2}e^{t_{2}}+2\alpha
_{1}\alpha _{2}e^{t_{2}}=0
\end{equation*}%
we get
\begin{equation*}
\alpha _{1}=\frac{e^{t_{1}}-2\alpha _{2}e^{t_{2}}+\alpha _{2}^{2}e^{t_{2}}}{%
2e^{t_{1}}-2\alpha _{2}e^{t_{2}}}.
\end{equation*}%
Similarly, from
\begin{equation*}
\frac{\partial \tilde{G}(\alpha )}{\partial \alpha _{2}}=\allowbreak
e^{t_{2}}\left( \alpha _{1}-1\right) \left( \alpha _{1}+2\alpha
_{2}-1\right) =0
\end{equation*}%
we get
\begin{equation*}
\alpha _{2}=\frac{1}{2}-\frac{1}{2}\alpha _{1}.
\end{equation*}%
Solving simultaneously we find
\begin{equation*}
\begin{aligned} \alpha _{1}^{\ast } & =\frac{1}{3e^{t_{2}}}\left(
-4e^{t_{1}}+3e^{t_{2}}+2\sqrt{e^{t_{1}}\left( 4e^{t_{1}}-3e^{t_{2}}\right)
}\right), \\ \alpha _{2}^{\ast } &= \frac{1}{3e^{t_{2}}}\left(
2e^{t_{1}}-\sqrt{e^{t_{1}}\left( 4e^{t_{1}}-3e^{t_{2}}\right) }\right).
\end{aligned}
\end{equation*}%
Note that for $e^{t_{1}}\rightarrow e^{t_{2}}$ we find
\begin{equation*}
\begin{aligned} \alpha _{1}^{\ast} &\rightarrow \frac{1}{3}, \\ \alpha
_{2}^{\ast} &=\frac{1}{2}-\frac{1}{2}\alpha _{1}^{\ast}\rightarrow
\frac{1}{3}. \end{aligned}
\end{equation*}%
Note also that if $\frac{3}{4}e^{t_{2}}<e^{t_{1}}<e^{t2}$ then $\alpha
_{1}^{\ast }<\alpha _{2}^{\ast }.$

II. From $\alpha _{1}=0,\alpha _{2}\neq 0$ and $\displaystyle \frac{%
\partial \tilde{G}(\alpha )}{\partial \alpha _{2}}=0$ we get
\begin{equation*}
\begin{aligned} \alpha _{1}^{\ast} &=0, \\ \alpha _{2}^{\ast} &=\frac{1}{2}.
\end{aligned}
\end{equation*}%
Note that constraining the growth rate as above is sufficient to render $%
\tilde{G}(\alpha )$ negative-semi definite. Alternatively one may add a
convex cost term $C(\alpha ).$ We have left aside these condition in what
follows to avoid notational clutter. By a straightforward computation, one
can easily check that the maximum of $\tilde{G}(\alpha )$ is attained for $%
\alpha ^{\ast }\in Int(K)$.

2. For $N=3$ we then have
\begin{equation*}
\begin{aligned} \tilde{G}(\alpha) = &\,y_{0}\left(\alpha _{1}(1-\alpha
_{1})e^{t_{1}}+\alpha _{2}(1-(\alpha _{1}+\alpha _{2}))e^{t_{2}}(1-\alpha
_{1}) \right.\\ &\left.+ \alpha _{3}(1-(\alpha _{1}+\alpha _{2}+\alpha
_{3}))e^{t_{3}}(1-\alpha_{1})(1-\alpha _{2})\right). \end{aligned}
\end{equation*}

I. From $\displaystyle \frac{\partial \tilde{G}(\alpha )}{\partial \alpha
_{1}}=0,$ we find
\begin{equation}
\alpha _{1}=\frac{e^{t_{1}}-2\alpha _{2}e^{t_{2}}-2\alpha
_{3}e^{t_{3}}+\alpha _{2}^{2}e^{t_{2}}+\alpha _{3}^{2}e^{t_{3}}-\alpha
_{2}\alpha _{3}^{2}e^{t_{3}}-\alpha _{2}^{2}\alpha _{3}e^{t_{3}}+3\alpha
_{2}\alpha _{3}e^{t_{3}}}{2e^{t_{1}}-2\alpha _{2}e^{t_{2}}-2\alpha
_{3}e^{t_{3}}+2\alpha _{2}\alpha _{3}e^{t_{3}}}.  \label{A}
\end{equation}%
From $\displaystyle \frac{\partial \tilde{G}(\mathbf{\alpha })}{\partial
\alpha _{2}}=0$ we find%
\begin{equation}
\alpha _{1}=\frac{1}{e^{t_{2}}-\alpha _{3}e^{t_{3}}}\left( e^{t_{2}}-2\alpha
_{2}e^{t_{2}}-2\alpha _{3}e^{t_{3}}+\alpha _{3}^{2}e^{t_{3}}+2\alpha
_{2}\alpha _{3}e^{t_{3}}\right) .  \label{B}
\end{equation}%
From $\displaystyle \frac{\partial \tilde{G}(\mathbf{\alpha })}{\partial
\alpha 3}=0$ we find%
\begin{equation}
\alpha _{1}=1-\alpha _{2}-2\alpha _{3}.  \label{C}
\end{equation}%
Combining (\ref{A}) and (\ref{B}) yields:%
\begin{eqnarray*}
&&\alpha _{2}=\frac{1}{6e^{t_{2}}-6\alpha _{3}e^{t_{3}}}\times \\
&&\left(
\begin{array}{c}
4t_{1}-\sqrt{%
\begin{array}{c}
16e^{2t_{1}}+\alpha _{3}^{2}e^{2t_{3}}+2\alpha _{3}^{3}e^{2t_{3}}+\alpha
_{3}^{4}e^{2t_{3}} \\
-12t_{1}e^{t_{2}}-16\alpha _{3}^{2}e^{t_{1}}e^{t_{3}}+12\alpha
_{3}^{2}e^{t_{2}}e^{t_{3}}-4\alpha _{3}e^{t_{1}}e^{t_{3}}%
\end{array}%
} \\
-5\alpha _{3}e^{t_{3}}+\alpha _{3}^{2}e^{t_{3}}%
\end{array}%
\right)
\end{eqnarray*}%
Combining (\ref{B}) and (\ref{C}) yields:
\begin{eqnarray*}
&&\alpha _{2}=\frac{1}{2e^{t_{2}}-2\alpha _{3}e^{t_{3}}}\times \\
&&\left(
\begin{array}{c}
2t_{1}-\sqrt{%
\begin{array}{c}
4e^{2t_{1}}+16e^{2t_{3}}e^{2t_{2}}+\alpha _{3}^{2}e^{2t_{3}}+6\alpha
_{3}^{3}e^{2t_{3}}+ \\
9\alpha _{3}^{4}e^{2t_{3}}-4e^{t_{1}}e^{t_{2}}-4\alpha
_{3}^{2}e^{t_{1}}e^{t_{3}}-4\alpha _{3}^{2}e^{t_{2}}e^{t_{3}}-24\alpha
_{3}^{3}e^{t_{2}}e^{t_{3}}%
\end{array}%
} \\
-4\alpha _{3}e^{t_{2}}-\alpha _{3}e^{t_{3}}+3\alpha _{3}^{2}e^{t_{3}}%
\end{array}%
\right)
\end{eqnarray*}%
equating the two we get the optimal $\alpha _{3}$.

Note that for $e^{t_{1}}\rightarrow e^{t_{2}}\rightarrow e^{t_{3}}$ we have
\begin{equation*}
\alpha _{2}=\frac{1}{6}\frac{-5\alpha _{3}+\alpha _{3}^{2}-\sqrt{-4\alpha
_{3}-3\alpha _{3}^{2}+2\alpha _{3}^{3}+\alpha _{3}^{4}+4}+4}{1-\alpha _{3}}
\end{equation*}%
and
\begin{equation*}
\alpha _{2}=\frac{1}{2}\frac{-5\alpha _{3}+3\alpha _{3}^{2}-3\sqrt{\alpha
_{3}^{2}-2\alpha _{3}^{3}+\alpha _{3}^{4}}+2}{1-\alpha _{3}}.
\end{equation*}%
This system is solved by
\begin{equation*}
\alpha _{i}^{\ast }=1/4,\quad i=1,2,3.
\end{equation*}%
II. a) From $\alpha _{1}=0,\alpha _{2}\neq 0,\alpha _{3}\neq 0$ and $%
\displaystyle \frac{\partial \tilde{G}(\alpha )}{\partial \alpha _{2}}=0$, $%
\displaystyle \frac{\partial \tilde{G}(\alpha )}{\partial \alpha _{3}}=0$ we
get
\begin{equation*}
\begin{aligned} &\alpha
_{3}=\frac{1}{3\tilde{t}_{3}}\left(2\tilde{t}_2-\sqrt{\tilde{t}_2\left(4%
\tilde{t}_2-3\tilde{t}_3\right)}\right),\\ &\alpha _{2}=1-2\alpha
_{3}=\frac{1}{3\tilde{t}_{3}}\left(3\tilde{t}_3-4\tilde{t}_2+2\sqrt{%
\tilde{t}_2\left(4\tilde{t}_2-3\tilde{t}_3\right)}\right). \end{aligned}
\end{equation*}%
Note that for $e^{t_{2}}\rightarrow e^{t_{3}}$
\begin{equation*}
\begin{aligned} \alpha _{1}^{\ast} &=0, \\ \alpha _{2}^{\ast} &\rightarrow
\frac{1}{3},\\ \alpha _{3}^{\ast} &\rightarrow \frac{1}{3}. \end{aligned}
\end{equation*}%
II. b) From $\alpha _{1}=0,\alpha _{2}=0,\alpha _{3}\neq 0$ and $%
\displaystyle \frac{\partial \tilde{G}(\alpha )}{\partial \alpha _{3}}=0$ we
get
\begin{equation*}
\begin{aligned} \alpha _{1}^{\ast} &=\alpha _{2}^{\ast}=0, \\ \alpha
_{3}^{\ast} &=\frac{1}{2}. \end{aligned}
\end{equation*}%
One can easily check that the maximum of $\tilde{G}(\alpha )$ is attained
for $\alpha ^{\ast }\in Int(K)$.
\qed

\paragraph{Proof of Lemma 7.}

Because payoffs are discounted there will be an
interior discount rate $\rho >0$ small enough, such that
\begin{eqnarray*}
&&\max_{\alpha _{i}\in \mathcal{A}}\int_{0}^{t_{1}}\left( e^{-\rho t}\alpha
_{_{i}}(t)\left( 1-\alpha _{i}(t)-\sum_{j\neq i}\alpha _{j}^{\ast
,r}(t)\right) f^{\alpha }(t-)-C(\alpha (t))\right) dt+ \\
&&\int_{t_{1}}^{T}\left( e^{-\rho t}\alpha _{i}(t)\left( 1-\alpha
_{i}(t)-\sum_{j\neq i}\alpha _{j}^{\ast }(t)\right) f^{\alpha }(t-)-C(\alpha
(t))\right) dt< \\
&&\int_{0}^{T}\left( e^{-\rho t}\alpha _{i}^{\ast ,r}(t)\left( 1-\alpha
_{i}^{\ast ,r}(t)-\sum_{j\neq i}\alpha _{j}^{\ast ,r}(t)\right) f^{\alpha
}(t-)-C(\alpha (t))\right) dt
\end{eqnarray*}%
$\forall i=\left \{ 1,...,I\right \} $ and deviation does not pay off with
grim punishments. Clearly also punishment that \textquotedblleft fit the
crime\textquotedblright \ as in Green and Porter (1984) may be employed.
\qed

\paragraph{Proof of Proposition 8.}

For the cartel behaviour the aggregate
statistic can be rewritten for each round $N>1$ as%
\begin{equation*}
\tilde{p}^{M,r}=1-\left( \frac{1}{I}\sum_{i=1}^{I}\alpha _{i}^{M,r}+\frac{1}{%
I}\sum_{i=1}^{I}\tilde{\varepsilon}_{i}\right)
\end{equation*}%
An individual deviation (w.l.o.g. by firm $i$)$~$changes the aggregate
statistic to%
\begin{equation*}
\tilde{p}^{d,r}=1-\left( \frac{1}{I}(\alpha _{i}^{d}+\frac{1}{I}%
\sum_{-i}^{I}\alpha _{-i}^{M,r})+\frac{1}{I}\sum_{i=1}^{I}\tilde{\varepsilon}%
_{i}\right)
\end{equation*}%
and so the visibility\ of individual (deviation) actions is decreasing in $%
\circ (1/I).$

Rewrite errors as%
\begin{equation*}
\frac{1}{I}\sum_{i=1}^{I}\tilde{\varepsilon}_{i}=\frac{\sqrt{Var(\varepsilon
)}}{\sqrt{I}}\frac{\sqrt{I}\left( \frac{1}{I}\sum_{i=1}^{I}\tilde{\varepsilon%
}_{i}\right) }{\sqrt{Var(\varepsilon )}}
\end{equation*}%
and note that by the Central Limit Theorem
\begin{equation*}
\frac{\sqrt{I}\left( \frac{1}{I}\sum_{i=1}^{I}\tilde{\varepsilon}_{i}\right)
}{\sqrt{Var(\varepsilon )}}\rightarrow N(0,1)
\end{equation*}%
so that the noise term decreases in $\circ (1/\sqrt{I})$ which prevents
individual detection in the limit. If it is individually optimal to go only
one round i.e. $N=1~$(e.g. if $C(\alpha (t))$ is not sufficiently convex),
then the limit principle holds \emph{irrespective} of the degree of noise.
\qed

\paragraph{Acknowledgements.}

The work by S. Ani\c{t}a and A.-M. Mo\c{s}neagu
was supported by the CNCS-UEFISCDI (Romanian National Authority for
Scientific Research) grant 68/2.09.2013, PN-II-ID-PCE-2012-4-0270:
\textquotedblleft Optimal Control and Stabilization of Nonlinear Parabolic
Systems with State Constraints. Applications in Life Sciences and
Economics\textquotedblright .

\pagebreak

\end{document}